\documentclass[12pt,leqno]{article}
\usepackage{amsmath,amsthm}

\def\init{\setcounter{equation}{0}}
\newtheorem{theorem}{Theorem}[section]

\newcommand{\R}{{\bf R}}

\newtheorem{lemma}{Lemma}[section]

\newcommand{\e}{{\varepsilon}}
\newcommand{\D}{\mathcal{D}}

\title{Inverse hyperbolic problems with 
time-dependent coefficients.
\author{G.Eskin, \ \ \  Department of Mathematics, UCLA,\\ Los Angeles,
CA 90095-1555, USA. \ E-mail: eskin@math.ucla.edu}
}

\begin{document}

\maketitle
\begin{abstract}
We consider the inverse problem for the second order 
self-adjoint hyperbolic equation
in a bounded domain in $\R^n$ with lower order terms depending analytically
on the time variable.   We prove that, assuming the BLR condition,
the time-dependent Dirichlet-to-Neumann
operator prescribed on a part of the boundary uniquely determines 
the coefficients of the hyperbolic equation up to a diffeomorphism and
a gauge transformation. As a by-product we prove a similar result for the
nonself-adjoint  hyperbolic operator with time-independent coefficients.
\end{abstract}

\section{Introduction.}
\label{section 1}
\init

Let $\Omega\subset \R^n$ be a bounded domain with a smooth boundary 
$\partial\Omega$ and let $\Gamma_0$ be an open subset of $\partial\Omega$.
Consider a hyperbolic equation in $\Omega\times(0,T_0)$ of the form:
\begin{eqnarray}                               \label{eq:1.1}
Lu\stackrel{def}{=}
\left(-i\frac{\partial}{\partial t}+A_0(x,t)\right)^2u(x,t)
\ \ \ \ \ \ \ \ \ \ \ \ \ \ \ \ \ \ \ \ \ \ 
\\
-\sum_{j,k=1}^n\frac{1}{\sqrt{g(x)}}
\left(-i\frac{\partial}{\partial x_j}+A_j(x,t)
\right)
 \sqrt{g(x)}g^{jk}(x)
\left(-i\frac{\partial}{\partial x_k}+A_k(x,t)\right)u
\nonumber
\\
-V(x,t)u=0,
\ \ \ \ \ \ \ \ \ \ \ \ 
\nonumber
\end{eqnarray}
where  
$ \|g^{jk}(x)\|^{-1}$ is the 
metric tensor
 in $\overline{\Omega}$,
$g(x)=\det\|g^{jk}\|^{-1}$ ,  $A_j(x,t),\ 0\leq j\leq n$,and $V(x,t)$ 
are smooth in $x\in\overline{\Omega}$ and real analytic in $t,\ t\in[0,T_0]$.
  We assume that
\begin{equation}                                 \label{eq:1.2}
u(x,0)=u_t(x,0)=0 \ \ \mbox{in}\ \ \ \Omega,
\end{equation}
and
\begin{equation}                                 \label{eq:1.3}
u\left|_{\partial\Omega\times(0,T_0)}\right. = f(x',t),
\end{equation}
where 
\begin{equation}                                 \label{eq:1.4}
\mbox{supp}\ f(x',t)\subset\overline{\Gamma_0}\times(0,T_0].
\end{equation}
Let $\Lambda f$ be the D-to-N (Dirichlet-to-Neumann) operator
on $\Gamma_0\times(0,T_0)$ :
\begin{equation}                                   \label{eq:1.5}
\Lambda f=\sum_{j,k=1}^n g^{jk}(x)\left(\frac{\partial u}{\partial x_j}
+iA_j(x,t)u\right)\nu_k\left
(\sum_{p,r=1}^ng^{pr}(x)\nu_p\nu_r\right)^{-\frac{1}{2}}
\left|_{\Gamma_0\times(0,T_0)}\right.,
\end{equation}
where $\nu=(\nu_1,...,\nu_n)$ is the unit exterior normal to 
$\partial\Omega$ with respect to the Euclidian metric.  We
shall study the inverse  
problem of the determination of the coefficients of (\ref{eq:1.1})
knowing the D-to-N operator on $\Gamma_0\times(0,T_0)$ for all smooth
$f$ with supports in $\overline{\Gamma_0}\times(0,T_0]$.

Let
\begin{equation}                               \label{eq:1.6}
y=y(x)
\end{equation}
be a diffeomorphism of $\overline{\Omega}$ onto some domain
$\overline{\Omega_0}$  such that $\partial\Omega_0\supset\Gamma_0$ and
\begin{equation}                               \label{eq:1.7}
y(x)=x \ \ \ \mbox{on}\ \ \ \Gamma_0.
\end{equation}
The equation (\ref{eq:1.1}) will have the following form in 
y-coordinates:
\begin{eqnarray}                               \label{eq:1.8}
L_0v\stackrel{def}{=}
\left(-i\frac{\partial}{\partial t}+A_0^{(0)}(y,t)\right)^2v(y,t)
\ \ \ \ \ \ \ \ \ \ \ \ \ \ \ \ \ \ \ \ \ \ 
\\
-\sum_{j,k=1}^n
\frac{1}{\sqrt{g_0(y)}}
\left(-i\frac{\partial}{\partial y_j}
+A_j^{(0)}(y,t)\right)
\sqrt{g_0(y)}
g_0^{jk}(y)
\left(-i\frac{\partial}{\partial y_k}
+   A_k^{(0)}(y,t)\right)v(y,t)
\nonumber
\\
-V^{(0)}(y)v(y,t)=0,
\ \ \ \ \ \ \ \ \ \ \ \ \ \ \ \ \ \
\nonumber 
\end{eqnarray}
where 
$v(y(x),t)=u(x,t)$,
\begin{equation}                                 \label{eq:1.9}
\|g_0^{jk}(y(x))\|=\left(\frac{\D y}{\D x}\right)\|g^{jk}(x)\|
\left(\frac{\D y}{\D x}\right)^T,
\end{equation}
\begin{equation}                                  \label{eq:1.10}
g_0(y)=\det \|g_0^{jk}(y)\|^{-1},
\end{equation}
$\frac{\D y}{\D x}$ is the Jacobi matrix of (\ref{eq:1.6}),
and
\begin{equation}                                   \label{eq:1.11}
A_0^{(0)}(y(x),t)=A_0(x,t),\ \ \ \ \ \ \ \ V^{(0)}(y(x),t)=V(x,t),
\end{equation}

\begin{equation}                                   \label{eq:1.12}
A_k(x,t)=\sum_{j=1}^nA_j^{(0)}(y(x),t)\frac{\partial y_j(x)}
{\partial x_k},\ \ \ 1\leq k\leq n. 
\end{equation}
Denote by $G_0(\overline{\Omega_0}\times [0,T_0])$ the group of
$C^\infty(\overline{\Omega_0}\times[0,T_0])$
complex valued functions $c(y,t)$ 
 such that $c(y,t)\neq 0$ and 
\begin{equation}                              \label{eq:1.13}  
c(y,t)=1\ \ \ \mbox{on}\ \ \Gamma_0\times[0,T_0].
\end{equation}

We shall call $G_0(\overline{\Omega_0}\times[0,T_0])$ the gauge group.
We say that the potentials $A_0^{(0)}(y,t),...,A_n^{(0)}(y,t),
V^{(0)}(y,t)$  and  $A_0^{(1)}(y,t),...,A_n^{(1)}(y,t),
V^{(1)}(y,t)$  are gauge equivalent if there exists 
$c(y,t)\in G_0(\overline{\Omega_0}\times[0,T_0])$
such that
\begin{equation}                                    \label{eq:1.14}
A_0^{(1)}(y,t)=A_0^{(0)}(y,t)-ic^{-1}(y,t)\frac{\partial c}{\partial t},
\ \ \ V^{(1)}(y,t)=V^{(0)}(y,t),
\end{equation}
\begin{equation}                                    \label{eq:1.15}
A_j^{(1)}(y,t)=A_j^{(0)}(y,t)-ic^{-1}(y,t)\frac{\partial c}{\partial y_j},
\ \ \ \ 1\leq j\leq n.
\end{equation}

Note that if
\begin{equation}                                     \label{eq:1.16}
v^{(1)}(y,t)=c^{-1}(y,t)v(y,t),
\end{equation}
then $v^{(1)}$ satisfies equation of the form (\ref{eq:1.8}) with
potentials $\{A_j^{(0)},V^{(0)}\}$ replaced by the gauge equivalent
$\{A_j^{(1)},V^{(1)}\},\ 0\leq j\leq n$.  We shall call 
the transformation (\ref{eq:1.16}) the gauge transformation.

Denote
 \begin{equation}                      \label{eq:1.17}
T_*=\sup_{x\in\overline{\Omega}}d(x,\Gamma_0),
\end{equation}
 where 
$d(x,\Gamma_0)$ is the distance in $\overline{\Omega}$ with respect to metric
tensor
$\|g^{jk}\|^{-1}$ from $x\in\overline{\Omega}$ to $\Gamma_0$.

We shall assume that the BLR condition (see [BLR]) is satisfied for
$t=T_{**}$.  This means roughly speaking that any bicharacteristics of
$L$ in $T_0^*(\Omega\times[0,T_{**}])$ intersects
$((\Gamma_0\times[0,T_{**}])\times\R^{n+1}\setminus\{0\}$.

Note that since BLR condition is determined by the geometry of $\Omega$ and 
$\Gamma_0$  and the second order terms of $L$ it holds when $[0,T_{**}]$
is replaced by $[t_0,T_{**}+t_0],\ \forall t_0>0$.

We shall  prove the following theorem:
\begin{theorem}                                    \label{theo:1.1}
Let $L$ and $L_0$ be two operators of the form (\ref{eq:1.1}) and
(\ref{eq:1.8}) 
in $\Omega$ and $\Omega_0$ respectively,
with coefficients analytic in $t\in [0,T_0]$  and let
$L$ and $L_0$ be formally self-adjoint, i.e coefficients $A_j,V$ and $A_j^{(0)},V^{(0)}$    
are real-valued, $0\leq j\leq n$.
   Suppose that
the BLR condition for $L$ is satisfied when $t=T_{**}$,  and
the D-to-N operators $\Lambda$ and $\Lambda^{(0)}$ 
coresponding to $L$ and $L^{(0)}$ respectively
are equal on
$\Gamma_0\times(0,T_0)$ for all $f$ satisfying (\ref{eq:1.4})
where $\Gamma_0\subset\partial\Omega\cap\partial\Omega_0$.
Suppose $T_0>2T_*+T_{**}$.  Then there exists a diffeomorphism $y=y(x)$
of $\overline{\Omega}$ onto $\overline{\Omega_0},\ y(x)=x$ on
$\Gamma_0$,  such that (\ref{eq:1.9}) holds.  Moreover,   
there exists a gauge transformation 
$c_0(x,t)\in G_0(\overline{\Omega}\times[0,T_0])$ 
such that
\begin{equation}                                    \label{eq:1.18}
c_0\circ y^{-1}\circ L_0=L
\end{equation}
in $\Omega\times(0,T_0)$.
\end{theorem}

{\bf Remark 1.1}
Let $L^*$ be formally adjoint operator to $L$. To prove
Theorem \ref{theo:1.1} we need to know in the addition to
D-to-N operator $\Lambda$ the D-to-N operator $\Lambda_*$
corresponding to $L^*$.  In the case when $L^*=L$ we have,
obviously,  that $\Lambda_*=\Lambda$.  When $A_0=0$ and 
the coefficients of $L$ are independent of $t$ one can show
( see,  for example, [KL1], [E1] ) that we can recover $\Lambda_*$ from
$\Lambda$.  
Therefore the proof and the result of Theorem \ref{theo:1.1} hold for the
case when $A_0=0$ and $A_j(x),\ 1\leq j\leq n,\ V(x)$ are complex-valued 
and independent of $t$, and this gives a new proof of the corresponding
result in [KL1].
\qed

The first inverse problem with boundary data on a part of the boundary
was considered in [I].  
The most general results were obtained by the BC-methods in the case
of self-adjoint hyperbolic operators with time-independent coefficients
( see [B1], [B2], [KKL], [KK] ).  The nonself-adjoint case with 
time-independent coefficient was considered in [B3], [KL1], [KL2], [KL3].
  The inverse problems for the wave equations  with time-dependent potentials were
considered  in [St],  [RS] (see also [I]).

The present paper is a generalization of the paper [E1].  We shall
widely use the notations,  results and the proofs from [E1].
Note that the case of equations with Yang-Mills potentials was
considered in [E2].

In \S 2 we recover the coefficients of $L$ (modulo a diffeomorphism 
and a gauge transformation) locally near $\Gamma_0$ (Theorem \ref{theo:2.1}).
Following [E1] in \S 3 we prove the global result.  In \S 4 we prove 
some lemmas used in \S 2.

\section{The local result.}
\label{section 2}
\init

Let $\Gamma$ be an open subset of $\Gamma_0$ and $U_0\subset \R^n$
be a neighborhood of $\Gamma$.  Let $(x',x_n)$  be coordinates in 
$U_0$ such that the equation of $\Gamma$ is $x_n=0$ and 
$x'=(x_1,...,x_{n-1})$
are coordinates on $\partial\Omega\cap U_0$.
Introduce semi-geodesic coordinates $y=(y_1,...,y_n)$ for $L$
in $U_0$  (c.f. [E1]):
\begin{equation}                       \label{eq:2.1}
y=\varphi(x),
\end{equation}
where $\varphi(x',0)=x'$,
\begin{equation}                       \label{eq:2.2}
\hat{g}^{nn}(y)=1,\ \ \hat{g}^{nj}=\hat{g}^{jn}=0,\ \ 1\leq j\leq n-1,
\end{equation}
$\|\hat{g}^{jk}(y)\|^{-1}$ is the metric tensor in the semi-geodesic
coordinates.  The equation (\ref{eq:1.1}) has the following form in
the semi-geodeic coordinates
\begin{eqnarray}                             \label{eq:2.3}
\hat{L}\hat{u}\stackrel{def}{=}\left(-i \frac{\partial}{\partial t}
+\hat{A}_0(y,t)\right)^2
\ \ \ \ \ \ \ \ \ \ \ \ \ \ \ \ \ \ \ \ \ \ \ \ \ \ \ \ \ \ \
\\
-\sum_{j,k=1}^n\frac{1}{\sqrt{\hat{g}(y)}}\left(-i\frac{\partial}{\partial y_j}
+\hat{A}_j(y,t)\right)
\sqrt{\hat{g}(y)}\hat{g}^{jk}(y)
\left(-i\frac{\partial}{\partial y_k}+\hat{A}_k(y,t)\right)\hat{u}(y,t)
\nonumber
\\
-\hat{V}(y,t)\hat{u}(y,t)=0,
\ \ \ \ \ \ \ \ \ \ \ \ \ \ \ \ \ \ \ \ \ \ \ \ \ \ \ \ \ \ \  \
\nonumber
\end{eqnarray}
where, as in \S 1, $\hat{u}(\varphi(x),t)=u(x,t),\ 
\hat{g}(y)=\det \|\hat{g}^{jk}(y)\|^{-1}$,
\begin{equation}                             \label{eq:2.4}
\hat{A}_0(\varphi(x),t)=A_0(x,t),\ \ \hat{V}(\varphi(x),t)=V(x,t),
\end{equation}
\begin{equation}                             \label{eq:2.5}
A_k(x,t)=\sum_{j=1}^n\hat{A}_j(\varphi(x),t)
\frac{\partial\varphi_j(x)}{\partial x_k} ,
\end{equation}
$\varphi(x)=(\varphi_1(x),...,\varphi_n(x))$.
As in [E1] let 
\begin{equation}                              \label{eq:2.6}
A_j'(y)=-\frac{i}{2}(\sqrt{\hat{g}})^{-1}\frac{\partial\sqrt{\hat{g}}}
{\partial y_j}
=-\frac{i\hat{g}_{y_j}}{4\hat{g}},\ \ 1\leq j\leq n.
\end{equation}
Then (\ref{eq:2.3}) can be rewritten in the form:
\begin{eqnarray}                                \label{eq:2.7}
\hat{L}\hat{u}=\left(- i\frac{\partial}{\partial t}
+\hat{A}_0(y,t)\right)^2
-
\left(-i\frac{\partial}{\partial y_n}+\hat{A}_n(y,t)
+A_n'(y) \right)^2\hat{u}
\\
-\sum_{j,k=1}^{n-1}\left(-i\frac{\partial}{\partial y_j}+\hat{A}_j
+A_j'(y)\right)
\hat{g}^{jk}\left(-i\frac{\partial}{\partial y_k} +
\hat{A}_k+A_k'\right)\hat{u}
-\hat{V}_1(y,t)\hat{u}=0,
\nonumber
\end{eqnarray}
where 
\begin{equation}                               \label{eq:2.8}
\hat{V}_1(y,t)=(A_n')^2+i\frac{\partial A_n'}{\partial y_n}
+\sum_{j,k=1}^{n-1}\left(\hat{g}^{jk}A_j'A_k'+
i\frac{\partial}{\partial y_j}(\hat{g}^{jk}A_k')\right)
+\hat{V}(y,t).
\end{equation}
Denote by $\hat{\Lambda}$ the D-to-N operator corresponding to
$\hat{L}$. We have
\begin{equation}                                   \label{eq:2.9}
\hat{\Lambda}f=
\left(\frac{\partial}{\partial y_n}
+i\hat{A}_n(y,t)
\right)\hat{u}(y,t)\left|_{y_n=0,0<t<T_0}\right. .
\end{equation}
Make the transformation $\hat{u}=(\hat{g}(y',y_n))^{-\frac{1}{4}}u'$.
Then $L'u'=0$,  where $L'$ is the same as (\ref{eq:2.7}) with 
$\hat{A}_j+A_j'$ replaced by $\hat{A}_j,\ 1\leq j\leq n$.
Denote by $\Lambda'$ the D-to-N operator 
$\Lambda'f'=\left(\frac{\partial}{\partial y_n}
+i\hat{A}_n(y,t)\right)u'\left|_{y_n=0,0<t<T_0}\right.$,  where
$L'u'=0,\ f'=u'\left|_{y_n=0,0<t<T_0}\right.$.
Since $\hat{\Lambda}$ determines $\hat{g}(y',0)$ and
$\frac{\partial\hat{g}(y',0)}{\partial y_n}$ (see Remark 2.2 in [E1])
we get that $\hat{\Lambda}$ determines $\Lambda'$ on 
$\Gamma\times(0,T_0)$ (c.f. [E1], (2.12), (2.13) ).

Analogously to (\ref{eq:1.14}), (\ref{eq:1.15}) 
$A_j^{(1)}(y,t),\ 0\leq j\leq n$, and $\hat{A}_j,
\ 0\leq j\leq n$,  are called gauge equivalent in 
$\overline{U_0}\times(0,T_0)$ if there exists 
$\hat{\psi}(y',y_n,t)\in C^\infty(\overline{U_0}\times[0,T_0]),$
\begin{equation}                                    \label{eq:2.10}
\hat{\psi}(y',0,t)=0,\ \ \ 0<t<T_0,
\end{equation}
such that
\begin{equation}                                   \label{eq:2.11}
\hat{A}_j^{(1)}(y,t)=\hat{A}_j(y,t)+
\frac{\partial\hat{\psi}}{\partial y_j},\ 1\leq j\leq n,\ 
\hat{A}_0^{(1)}=\hat{A}_0+\frac{\partial\hat{\psi}}{\partial t}.
\end{equation}
Here $c=e^{i\hat{\psi}}$.  We shall choose $\hat{\psi}(y,t)$ such that
\begin{equation}                                   \label{eq:2.12}
\hat{A}_0-\hat{A}_n(y,t)+\frac{\partial\hat{\psi}}{\partial t}
-\frac{\partial\hat{\psi}}{\partial y_n}=0,
\end{equation}
i.e.
\begin{equation}                                \label{eq:2.13}
\hat{A}_0^{(1)}=\hat{A}_n^{(1)} \ \ \ \mbox{in}\ \ \ 
\overline{U_0}\times[0,T_0].
\end{equation}
Note that  $\hat{\psi}(y',y_n,t)$ is analytic in $t$ since
$\hat{A}_j(y,t),\ 0\leq j\leq n,$  are analytic in $t$.
Therefore $\hat{A}_j^{(1)},\ 0\leq j\leq n,$  are also analytic in $t$.

Let
\begin{equation}                               \label{eq:2.14}
s=t-y_n,\ \ \tau=T-t-y_n.
\end{equation}
Note that
\begin{equation}                               \label{eq:2.15}
\hat{u}_s=\frac{1}{2}(\hat{u}_t-\hat{u}_{y_n}),\ \ \hat{u}_\tau=
-\frac{1}{2}(\hat{u}_t+\hat{u}_{y_n}).
\end{equation}
Substituting $u'=e^{i\hat{\psi}}u_1$ in $L'u'=0$
we get the following equation for $u_1(y,t)$:
\begin{eqnarray}                                \label{eq:2.16}
L_1u_1\stackrel{def}{=}\left(-i \frac{\partial}{\partial t}
+\hat{A}_n^{(1)}(y,t)\right)^2u_1
-
\left(-i\frac{\partial}{\partial y_n}+\hat{A}_n^{(1)}(y,t)
 \right)^2u_1
\\
-\sum_{j,k=1}^{n-1}\left(-i\frac{\partial}{\partial y_j}
+\hat{A}_j^{(1)}(y,t)
\right)
\hat{g}^{jk}\left(-i\frac{\partial}{\partial y_k} +
\hat{A}_k^{(1)}(y,t)\right)u_1
-\hat{V}_1(y,t)u_1=0,
\nonumber
\end{eqnarray}
In $(s,\tau,y')$ coordinates we have
\begin{eqnarray}                                \label{eq:2.17}
L_1u_1=
4u_{1s\tau}-4i\hat{A}_n^{(1)}u_{1s}
\ \ \ \ \ \ \ \ \ \ \ \ \ \ \ \ \ \ \ \ \ \ \ \ \ \ \ 
\\
-\sum_{j,k=1}^{n-1}\left(-i\frac{\partial}{\partial y_j}
+\hat{A}_j^{(1)}(y,t)
\right)
\hat{g}^{jk}\left(-i\frac{\partial}{\partial y_k} +
\hat{A}_k^{(1)}(y,t)\right)u_1
-\hat{V}_2(y,t)u_1=0,
\nonumber
\end{eqnarray}
where
\[
\hat{V}_2=\hat{V}_1+2i\hat{A}_{ns}^{(1)}.
\]
Note that the D-to-N operator corresponding to $L_1$ is
\begin{equation}                           \label{eq:2.18}
\Lambda^{(1)}(u_1\left|_{y_n=0,0<t<T_0}\right.)
 =\left(\frac{\partial u_1}{\partial y_n} + 
i\hat{A}_n^{(1)}u_1\right)\left|_{y_n=0,0< t<T_0}\right. ,
\end{equation}
and
\begin{equation}                          \label{eq:2.19}
u_1(y,t)=e^{-i\hat{\psi}(y,t)}
(\hat{g}(y',y_n))^{\frac{1}{4}}u(\varphi^{-1}(y),t),
\end{equation}
where $Lu=0$.

Let $\Delta_{1s_0}=\Gamma\times[s_0,T]$ where $T_1\leq s_0 < T$.  Denote by 
$D(\Delta_{1s_0})$ the forward domain of influence of 
$\Delta_{1s_o}$ in the half-space $y_n\geq 0$. 
 Let $\Gamma^{(2)}=\{y':(y',y_n,t)\in D(\Delta_{1T_1}), y_n=0,t=T\}$
and let $\Delta_{2s_0}=\Gamma^{(2)}\times[s_0,T]$.
Denote by $D(\Delta_{2s_0})$ the forward domain of influence of
$\Delta_{2s_0}$ for $y_n\geq 0$.  Let
$\Gamma^{(3)}=\{y':(y',y_n,t)\in D(\Delta_{2T_1}),y_n=0,t=T\}$ and let
$\Delta_{3s_0}=\Gamma^{(3)}\times[s_0,T]$.  Denote by $D(\Delta_{3s_0})$ the
forward domain  of influence of $\Delta_{3s_0}$ in the half-space
$y_n\geq 0$.
We assume that $T-T_1$  is small such that 
$D(\Delta_{3T_1})\subset \Gamma_0\times [T_1,T]$
for $t\leq T$,  the semigeodesic coordinates are defined  in 
$D(\Delta_{3T_1}),t\leq T$,  and $D(\Delta_{3T_1})\setminus \{y_n=0\}$
does not intersect  $\partial\Omega\times [T_1,T]$.  Denote by 
$Y_{js_0}$ the intersection of $D(\Delta_{js_0})$ with the plane
$T-t-y_n=0,1\leq j\leq 3$.   Let $X_{js_0}$ be  the part of $D(\Delta_{js_0})$
below $Y_{js_0}$ and let $Z_{js_0}=\partial X_{js_0}\setminus (Y_{js_0}\cup
\{y_n=0\}),1\leq j\leq 3$.

Suppose $L_1u_1=0$ for $y_n>0,\ t<T,\ u_1=u_{1t}=0$  for $t=T_1,\ y_n>0,\ 
u_1\left|_{y_n=0,T_1<t<T}\right.=f.$
Let $L_1^*$ be the operator formally adjoint to $L_1$.
Note that $L_1^*$ has the form (\ref{eq:2.16})  with $\hat{A}_j^{(1)},
0\leq j\leq n, \hat{V}_1(y,t)$ replaced by 
$\overline{\hat{A}_j^{(1)}},
0\leq j\leq n, \overline{\hat{V}_1(y,t)}$.
The D-to-N operator $\Lambda_*^{(1)}$ corresponding to $L_1^*$ has
the form 

\[
\Lambda_*^{(1)}g =\left(\frac{\partial v_1}{\partial y_n} + 
i\overline{\hat{A}_n^{(1)}}(y,t)v_1\right)\left|_{y_n=0,T_1< t<T}\right. .
\]
where $L_1^*v_1=0$ for $y_n>0, \ t<T,\ v_1\left|_{y_n=0,T_1<t<T}\right. =g,
v_1=v_{1t}=0$ when $t=T_1, y_n>0$.

We assume that $\Lambda_*^{(1)}$ can be determined if we know $\Lambda^{(1)}$.
This is obviously true when $L_1$ is formally self-adjoint.
Then $\Lambda_*^{(1)}=\Lambda^{(1)}$.
Note that $L_1$ is self-adjoint if $L$ is self-adjoint.
  Note that $\Lambda^{(1)}$  determines $\Lambda_*^{(1)}$ also when 
$\hat{A}_0^{(1)}=0,\ A_j^{(1)},\ 
1\leq j\leq n, \hat{V}_1$ are independent of $t$ (c.f. [KL1]).

Consider the identity
\begin{equation}                              \label{eq:2.20}
0=(L_1u_1,v_1)-(u_1,L_1^*v_1),
\end{equation}
where 
\begin{equation}                              \label{eq:2.21}
(u_1,v_1)=\int_{X_{3T_1}}u_1(y,t)\overline{v_1(y,t)}dydt.
\end{equation}
We assume that $\mbox{supp}\ f$ and $\mbox{supp}\ g$ are contained
in $\Delta_{3T_1}$.

Integrating by parts we get
\begin{eqnarray}                             \label{eq:2.22}
\int_{X_{3T_1}}\left(-i\frac{\partial }{\partial t} + 
\hat{A}_n^{(1)}\right)^2u_1\overline{v_1}dy'dy_ndt
\ \ \ \ \ \ \ \ \ \ \ \ \ \ \ \ \ \ \ \ \ \ \ \ \ \ \ \ \ \ 
\\
=\int_{X_{3T_1}}\left(-i\frac{\partial }{\partial t} + \hat{A}_n^{(1)}\right)u_1
\overline{\left(-i\frac{\partial}{\partial t}+
\overline{\hat{A}_n^{(1)}}\right)
}v_1dydt
-i\int_{Y_{3T_1}}\left(-i\frac{\partial }{\partial t} + 
\hat{A}_n^{(1)}\right)u_1\overline{v_1} dy'dy_n.
\nonumber
\end{eqnarray}
We used that $u_1$ and $v_1$ are zero on $Z_{3T_1}$.
Analogously
\begin{eqnarray}                             \label{eq:2.23}
\int_{X_{3T_1}}\left(-i\frac{\partial }{\partial y_n} + 
\hat{A}_n^{(1)}\right)^2u_1\overline{v_1}dy'dy_ndt
 \ \ \ \ \ \ \ \ \ \ \ \ \ \ \ \ \ \ 
\\
=\int_{X_{3T_1}}\left(-i\frac{\partial }{\partial y_n} + \hat{A}_n^{(1)}\right)u_1
\overline{\left(-i\frac{\partial}{\partial y_n}+
\overline{\hat{A}_n^{(1)}}\right)
}v_1dydt
\nonumber
\\
-i\int_{Y_{3T_1}}\left(-i\frac{\partial }{\partial y_n} + 
\hat{A}_n^{(1)}\right)u_1\overline{v_1} dy'dt.
+i\int_{\Delta_{3T_1}}\left(-i\frac{\partial }{\partial y_n} + 
\hat{A}_n^{(1)}\right)u_1\overline{v_1} dy'dt.
\nonumber
\end{eqnarray}
Analogously,  integrating by parts other terms in (\ref{eq:2.20}) and
taking into account that $u_1=v_1=0$ on $Z_{3T_1}$  we get
\begin{eqnarray}                              \label{eq:2.24}
0=(L_1u_1,v_1)-(u_1,L_1^*v_1)=
-\int_{Y_{3T_1}}(u_{1s}\overline{v_1}-u_1\overline{v_{1s}})dy'ds
\\
+\int_{\Delta_{3T_1}}
\left[\left(\frac{\partial }{\partial y_n} + 
i\hat{A}_n^{(1)}\right)u_1\overline{v_1}
-u_1\overline{\left(\frac{\partial }{\partial y_n} + 
i\overline{\hat{A}_n^{(1)}}\right)v_1}\ \ \right]dy'dt.
\nonumber
\end{eqnarray}
We used in (\ref{eq:2.24}) the change of variables (\ref{eq:2.14}) and 
(\ref{eq:2.15}).  Note that the integrals over $Y_{3T_1}$ containing
$\hat{A}_n^{(1)}$ are cancelled.  The second integral in (\ref{eq:2.24})
has the form
\begin{equation}                               \label{eq:2.25}
(\Lambda^{(1)}f,g)-(f,\Lambda_*^{(1)}g),
\end{equation}
where $u_1=f$ and $v_1=g$ on $\Delta_{3T_1}$.
Denote
\begin{equation}                               \label{eq:2.26}
A_0(u_1,v_1)= 
\int_{Y_{3T_1}}(u_{1s}\overline{v_1}-u_1\overline{v_{1s}})dy'ds.
\end{equation}
It follows from (\ref{eq:2.24}),  (\ref{eq:2.25})  that $A_0(u_1,v_1)$
is determined by the D-to-N  operator.  Integrating by parts in
(\ref{eq:2.26}) we get
\[
A_0(u_1,v_1)=2\int_{Y_{3T_1}}u_{1s}\overline{v}dsdy' -
\int_{\gamma_{3T_1}}u_1(y',0,T)
\overline{v_1(y',0,T)}dy',
\]
where $\gamma_{js_0}=\partial Y_{js_0}\cap\{y_n=0\},\ 1\leq j\leq 3$.
Since $u_1(y',0,T)=f(y',T),\ v_1(y',0,T)=g(y',T)$  we get that
\begin{equation}                              \label{eq:2.27}
A(u_1,v_1)\stackrel{def}{=}2\int_{Y_{3T_1}}u_{1s}\overline{v_1}dy'ds
\end{equation}
is determined by the D-to-N operator on $\Delta_{3T_1}$.

Denote by $\stackrel{\circ}{H^1}(\Delta_{js_0})$ the subspace of the 
Sobolev space $H^1(\Delta_{js_0})$
consisting of functions equal to zero on $\partial\Delta_{js_0}$
and by $H_0^1(\Delta_{js_0})$ the subspace of $H^1(\Delta_{js_0})$
consisting of functions equal to zero on
$\partial\Delta_{js_0}\setminus\{t=T\}, j=1,2,3$.    Also
denote by $H_0^1(Y_{js_0})$ the subspace of $H^1(Y_{js_0})$
consisting of functions equal to sero on $\partial Y_{js_0}\setminus
\{y_n=0\}$ and by
$\stackrel{\circ}{H^1}(Y_{1s_0})$ the subspace of $H^1(Y_{1s_0})$ of
 functions equal to zero on 
$\partial Y_{js_0}$.  Here
$s_0\in[T_1,T)$.

For the convenience we shall often denote by $u^f$ (correspondently
$v^g$) the solutions of $L_1u^f=0,\ u^f\left|_{y_n=0}\right.=f,\ 
u^f=u_t^f=0$ when $t=T_1$ (correspondently $L_1^*v^g=0,
\ v^g\left|_{y_n=0}=g,\right.
\ v^g=v_t^g=0$ when $t=T_1$).

Denote by $R_{js_0}$ the following subset of $Y_{js_0}:R_{js_0}=
\Gamma^{(j)}\times[s_0,T],1\leq j\leq 3$,  where $\Gamma^{(1)}=\Gamma$.
Note that $(Y_{1T_1}\cap\{s\geq s_0\})\subset R_{2s_0}\subset Y_{2s_0}
\subset R_{3s_0}\subset Y_{3s_0}$  for any $s_0\in [T_1,T)$.
\begin{lemma}                          \label{lma:2.1}
For any smooth $f\in H_0^1(\Delta_{1T_1})$ and any 
$s_0\in [T_1,T)$ there exists $u_0\in 
H_0^1(R_{2s_0})$ such that
\begin{equation}                              \label{eq:2.28}
A(u^f,v')=A(u_0,v')
\end{equation}
for all $v'\in H_0^1(Y_{3s_0})$.
\end{lemma}

{\bf Proof:}
It follows from (\ref{eq:2.28}) that $u_s^f=u_{0s}$ in $Y_{3s_0}$.
Let $w_1$ be such that $w_{1s}=0$ in $R_{2s_0},\ w_1
=u^f $ when $s=s_0, y'\in\Gamma^{(2)}$.  Then
$u_0=u^f-w_1\in H_0^1(R_{2s_0})$  is the unique solution of (\ref{eq:2.28}).
\qed

Extend $u_0$ by zero in $Y_{3s_0}\setminus R_{2s_0}$. 
Then $u_0\in H_0^1(Y_{3s_0})$.
 
\begin{lemma}                               \label{lma:2.2}
For any $v_j\in H_0^1(R_{js_0})$ there exists a sequence $u^{f_{nj}}\in
H_0^1(Y_{js_0})$ that converges to $v_j$  in 
$H_0^1(Y_{js_0})$.  Here $j=1,2,3.$ 
\end{lemma}
In the case when coefficients of $L$ do not depend on $t$ this lemma
was proven in [E1].  Note that the proof in [E1] works also for
nonself-adjoint $L$.  In the case of the time-dependent coefficients
Lemma \ref{lma:2.2} will be proven in \S 4.
\qed

Assume that BLR condition is satisfied for 
$L^{(1)}$  when
$t=T_{**}$.  Then for any
$T>T_{**}$ the map of $f\in H_0^1(\Gamma_0\times (0,T))$ to
$\{u^f(x,T),u_t^f(x,T)\}\in H^1(\Omega)\times L_2(\Omega)$ is onto
(see [BLR]).  We always assume the zero initial conditions when
$t=0,\ x\in \Omega$.

Note that (c.f. [H])
\[
\|u^f(\cdot,T)\|_{1,\Omega}+\|u_t^f(\cdot,T)\|_{0,\Omega}^2
\leq C\|f\|_{1,\Gamma_0\times(0,T)}.
\]
Therefore by the closed graph theorem we have 
that
\begin{equation}                              \label{eq:2.29}
\inf_{\mathcal{F}}\|f'\|_{1,\Gamma_0\times(0,T)}\leq 
C(\|u^f(\cdot,T)\|_{1,\Omega}^2+\|u_t^f(\cdot,T)\|_{0,\Omega}^2),
\end{equation}
where $\mathcal{F}\subset H_0^1(\Gamma_0\times(0,T))$ consists of all
$f'$ such that $u^{f'}(x,T)=u^f(x,T),\\
u_t^{f'}(x,T)=u_t^f(x,T),\ x\in\Omega$.

Let $L^{(i)},\ i=1,2,$  be two operators having the same D-to-N operator
on $\Gamma_0\times(0,T), \ T>T_{**}$.  Let $X_{js_0}^{(i)},\ Y_{js_0}^{(i)}, j=1,2,3$
be corresponding to $L_1^{(i)},\ i=1,2.$

Since $\Lambda^{(1)}=\Lambda^{(2)}$ on $\Delta_{3s_0}$ we get that
$D^{(1)}(\Delta_{1s_0})\cap\{y_n=0\}=D^{(2)}(\Delta_{1s_0})\cap\{y_n=0\}$ 
where $D^{(i)}(\Delta_{1s_0})$ is the forward domain of influence of 
$L^{(i)}$ in the half-space $y_n\geq 0,i=1,2$. (see [E1], Lemma 2.4).  Therefore 
$\Gamma_1^{(2)}=\Gamma_2^{(2)}$,  i.e. $\Delta_{2s_0}^{(1)}=\Delta_{2s_0}^{(2)}$.
Analogously one proves that $\Delta_{3s_0}^{(1)}=\Delta_{3s_0}^{(2)}$.
Therefore $\Delta_{js_0}^{(1)}=\Delta_{js_0}^{(2)}=\Delta_{js_0},j=1,2,3.$

\begin{lemma}                            \label{lma:2.3}
Assume that $L^{(i)},\ i=1,2,$
satisfy the BLR condition for $T>T_{**}$.  Let $u_i^f,\ i=1,2,$
be the solutions of $L_1^{(i)}u_i^f=0$ in $X_{2s_0}^{(i)}$ such that 
$u_1^f=u_2^f=f$ on 
$\Delta_{2s_0},\ \mbox{supp\ }f\subset \Delta_{2s_0},\ 
f\in H_0^1(\Delta_{2s_0}),\ s_0\in [T_1,T)$.
Then there exists constants $C_1$ and $C_2$ such that
\begin{equation}                        \label{eq:2.30}
C_1\|u_1^f\|_{1,Y_{2s_0}^{(1)}}\leq \|u_2^f\|_{1,Y_{2s_0}^{(2)}}\leq
C_2\|u_1^f\|_{1,Y_{2s_0}^{(1)}}
\end{equation}
\end{lemma}

{\bf Proof:} 
Let $\mathcal{F}$  be the same as above. 
Note that $\mbox{supp\ }u_1^f(x,T)$   
and $u_{1t}^f(x,T)$ are contained in
$ D(\Delta_{2s_0})\cap
\{t=T\}$.
Let $\Gamma_2,\Gamma_3,\Gamma_4$ be the same as in \S 4.

Note that $D(\Delta_{2s_0})\cap\{t=T\}\subset\Gamma_4$ and 
$Y_{2s_0}^{(1)}\subset\Gamma_2$.
The following estimate will be proven in \S 4 (c.f. Lemma 3.1 in [E1]):
\begin{eqnarray}                         \label{eq:2.31}
C_1(\|u_1^f\|_{1,\Gamma_2}^2+\|u_1^f\|_{1,\Gamma_3}^2)
\ \ \ \ \ \ \ \ \ \ \ \ \ \ \ \ \ \ \ \ \ 
\\
\leq \|u_1^f\|_{1,\Gamma_4}^2  + \|u_{1t}^f\|_{0,\Gamma_4}^2
\leq C_2(\|u_1^f\|_{1,\Gamma_2}+\|u_1^f\|_{1,\Gamma_3}^2),
\nonumber
\end{eqnarray}
assuming that $u_1^f$ has a compact support in $y'$.  Note that 
since $\mbox{supp\ }f\subset \Delta_{2s_0}$ we have that 
$u_1^f=0$ on $\Gamma_3$. Therefore (\ref{eq:2.31}) implies
\begin{equation}                       \label{eq:2.32}
\|u_1^f\|_{1,\Gamma_4}^2 + \|u_{1t}^f\|_{0,\Gamma_4}^2
\leq C\|u_1^f\|_{1,Y_{2s_0}^{(1)}} .
\end{equation}
For any $f'\in \mathcal{F}$ we have $u_1^{f'}=u_1^f,\ u_{1t}^{f'}=u_{1t}^f$
on $\Gamma_4$.  Therefore $u_1^{f'}|_{Y_{2s_0}^{(1)}}=u_1^{f}|_{Y_{2s_0}^{(1)}}$
by the domain of dependence arguments.

It follows from (\ref{eq:2.29}) that there exists $f_0\in \mathcal{F}$
such that
\begin{equation}                           \label{eq:2.33}
\|f_0\|_{1,\Gamma_0\times(0,T)}\leq C_1(\|u_1^f\|_{1,\Gamma_4}+
\|u_{1t}^f\|_{0,\Gamma_4}).
\end{equation}
Note that the solutions $u^f$ of $L^{(1)}u^f=0$ and $u_1^f$ of $L_1^{(1)}u_1^f=0$
in $D(\Delta_{2s_0})$ are related by (\ref{eq:2.19}).

Now we shall show that $u_2^{f'}|_{Y_{2s_0}^{(2)}}=u_2^{f}|_{Y_{2s_0}^{(2)}}$
for all $f'\in \mathcal{F}$.  
Note that the Green's formula (\ref{eq:2.24}) holds for any $u_1^{f'},\ v_1^g$
where $f'\in H_0^1(\Gamma_0\times(0,T_0)),\ g\in H_0^1(\Delta_{3s_0}),\ 
L_1u_1^{f'}=0,\ L_1^*v_1^g=0$ 
in $D(\Delta_{3s_0}),\ L_1$ is the same as in (\ref{eq:2.16}),  $Lu^{f'}=0$
in $\Omega\times(0,T_0),\ u^{f'}=u_t^{f'}=0$ for $t=0$,       $u^{f'}$
and $u_1^{f'}$ are related by (\ref{eq:2.19})  in $D(\Delta_{3s_0})$.
Since $\Lambda^{(1)}=\Lambda^{(2)}$ on
$\Gamma_0\times(0,T_0)$  we get from 
 (\ref{eq:2.19}) that  (\ref{eq:2.25}) is the same for
$i=1,2$.  Therefore 
\[
(u_{1s}^{f'},v_1^g)=(u_{2s}^{f'},v_2^g)
\]
for all $f'\in\mathcal{F}$ and all $g\in H_0^1(\Delta_{3s_0})$.
Since $u_1^{f'}|_{Y_{2s_0}^{(1)}}=u_1^{f}|_{Y_{2s_0}^{(1)}}$
for all $f'\in\mathcal{F}$ we get that
\[
(u_{2s}^{f'},v_2^g)=(u_{2s}^f,v_2^g)
\]
for all $g\in H_0^1(\Delta_{3s_0})$.
By Lemma \ref{lma:2.2}  with $j=3$  
$\{v_2^g\}$ are dense
in $H_0^1(R_{3s_0}^{(2)})$.  Therefore we have that $u_{2s}^{f'}=u_{2s}^f$
in $Y_{2s_0}^{(2)}$ for any $f'\in\mathcal{F}$,  since $R_{3s_0}^{(2)}
\supset Y_{2s_0}^{(2)}$.  Note that 
$u_2^{f'}\left|_{t=T}\right.
=f'(y',T)$ and $u_1^{f'}(y',0,T)=f'(y',T)$.
Therefore $u_2^{f'}|_{t=T}=u_2^f|_{t=T}$ since 
$u_1^f(y',0,T)=u_1^{f'}(y',0,T)$.  Therefore 
$u_2^{f'}|_{Y_{2s_0}^{(2)}}=u_2^{f}|_{Y_{2s_0}^{(2)}}$.

Estimate (\ref{eq:2.31}) is proven in \S 4 under the assumption that 
$u_2^{f_0}(y,T),\ u_{2t}^{f_0}(y,T)$ have  compact supports in $y'$.
Denote by $D_1$ the intersection of forward domain of influence of
$Y_{2s_0}^{(2)}$ with the plane $t=T,\ y_n \geq 0$.  Note that
$D_1\subset \Gamma_4$.  Let $\hat{u}_2(y,T),\ \hat{u}_{2t}(y,T)$
be  extensions of $u_2^{f_0}(y,T),\ u_{2t}^{f_0}(y,T)$
from $D_1$ to a neighborhood $\hat{D}_1$ of $D_1$ in $\Gamma_4$ such that
\[
\|\hat{u}_2\|_{1,\hat{D}_1}^2 + 
\|\hat{u}_{2t}\|_{0,\hat{D}_1}^2
\leq 2(\|u_2^{f_0}\|_{1,D_1}^2 + \|u_{2t}^{f_0}\|_{0,D_1}^2),
\]
and $\hat{u}_2, \ \hat{u}_{2t}$ are zero outside of $\hat{D}_1$.
Applying Lemma \ref{lma:4.1} to $\hat{u}_2,
\hat{u}_{2t}$ we get that
\begin{eqnarray}                          \label{eq:2.34}
\|u_2^{f_0}\|_{1,Y_{2s_0}^{(2)}}\leq 
C(\|\hat{u}_2\|_{1,\hat{D}_0}^2+\|\hat{u}_{2t}\|_{0,\hat{D}_0})
\\
\leq 
2C(\|u_2^{f_0}\|_{1,D_0}^2+\|u_{2t}^{f_0}\|_{0,D_0})
\nonumber
\\
\leq C_1(\|v_2^{f_0}\|_{1,\Omega}^2+\|v_{2t}^{f_0}\|_{0,\Omega})
\leq C_2\|f_0\|_{1,\Gamma_0\times(0,T)}^2,
\nonumber
\end{eqnarray}
where $v_2^{f_0}$ is the solution of $L^{(2)}v_2^{f_0}=0$ in 
$\Omega^{(2)}\times(0,T_0)$ corresponding to
$f_0\in H_0^1(\Gamma_0\times(0,T_0)),\ u_2^{f_0}$  is the solution of
$L_1^{(2)}u_2^{f_0}=0$ in $D(\Delta_{2s_0}),\ u_2^{f_0}
=e^{-i\psi^{(2)}}(\hat{g}_2)^{\frac{1}{4}}\varphi_2^{-1}\circ v_2^g$
(see (\ref{eq:2.19}) ).
We used
in (\ref{eq:2.34}) that the domain of dependence of $Y_{2s_0}^{(2)}$ 
intersected with $t=T$
is 
$D_1$ and therefore $u_2^{f_0}|_{Y_{2s_0}^{(2)}}$ does not depend on
the extensions $\hat{u}_2,\ \hat{u}_{2t}$.
Combining (\ref{eq:2.32}), (\ref{eq:2.33}) and (\ref{eq:2.34}) we  get
the right half of the inequality (\ref{eq:2.30}).

The inequality $\|u_1^f\|_{1,Y_{2s_0}^{(1)}}\leq C\|u_2^f\|_{1,Y_{2s_0}^{(2)}}$
also holds if we assume that $L^{(2)}$ satisfies the BLR property
too.
However,  to prove Theorem \ref{theo:1.1} we need only BLR condition for 
$L^{(1)}$.

\begin{lemma}                                \label{lma:2.4}
Suppose 
$g\in H_0^1(\Delta_{1T_1})$ is smooth.
Let $u_0$ be the same as in Lemma \ref{lma:2.1}.  Then the inner
product $(u_{0s},v^g)$ is uniquely determined by the D-to-N operator.
\end{lemma}
Note that $(u_{1s}^f,v_1^g)=(u_{2s}^f,v_2^g)$
are determined by D-to-N operator (see (\ref{eq:2.24}) where
$u_i^f,\ v_i^g$ correspond to $L_1^{(i)},\ (L_1^{(i)})^*,\ i=1,2$).  
Since $\{u_1^f\}, f\in H_0^1(\Delta_{2s_0})$ are dense in 
$H_0^1(R_{2s_0}^{(1)})$ there exists a sequence  $u_1^{f_n}$ convergent to
$u_0^{(1)}$ in $H_0^1(Y_{2s_0}^{(1)})$.
Here $u_0^{(i)},\ i=1,2,$  are the same as in (\ref{eq:2.28}) for 
$i=1,2,\ f_n\in H_0^1(\Delta_{2s_0})$.
By Lemma \ref{lma:2.3} $\{u_2^{f_n}\}$ converges in $H_0^1(Y_{2s_0}^{(2)})$
to some function $w_2\in H_0^1(Y_{2s_0}^{(2)})$.  Passing to the limit
when $n\rightarrow \infty$ we get 
\begin{equation}                         \label{eq:2.35}
(u_{0s}^{(1)},v_1^g)=(w_{2s},v_2^g).
\end{equation}

It follows from (\ref{eq:2.28}) that 
$(u_{0s}^{(1)},v_1^g)=(u_{0s}^{(2)},v_2^g)$ for any 
$g\in H_0^1(\Delta_{3s_0})$.
Compairing with (\ref{eq:2.35})
we get $(u_{0s}^{(2)},v_2^g)=(w_{2s}^{(2)},v_2^g)$.  Since, 
by Lemma \ref{lma:2.2} with $j=3$, $\{v_2^g\}$  are dense in 
 $H_0^1(R_{3s_0}^{(2)})\supset 
H_0^1(Y_{2s_0}^{(2)})$ we have 
$w_{2s}=u_{0s}^{(2)}$.  
Since $w_2=u_0^{(2)}=0$ on $\partial Y_{2s_0}^{(2)}\setminus\{t=T\}$
we get $u_0^{(2)}=w_2$ in $Y_{2s_0}^{(2)}$.
Therefore $(u_{0s}^{(1)},v_1^g)=(u_{0s}^{(2)},v_2^g)$ for all 
$g\in H_0^1(\Delta_{1T_1})$ since 
$Y_{1T_1}^{(i)}\cap\{s\geq s_0\}\subset R_{2s_0}^{(i)}\subset
Y_{2s_0}^{(i)},i=1,2, $
  i.e. 
 $A(u_0,v^g)=2(u_{0s},v^g)$ is uniquely determined
by the D-to-N. 
\qed

Therefore 
$A_1(u^f,v^g)\stackrel{def}{=}A(u^f,v^g)-A(u_0,v^g)$
is determined by the D-to-N operator  when
$g$ and $f$ are smooth and belong to
$ \stackrel{\circ}{H^1}(\Delta_{1T_1})$.

Since $u_s^f-u_{0s}=0$ in $R_{2s_0}$ and $u_0= 0$ in 
for $s\leq s_0$ we have
\begin{equation}                         \label{eq:2.36}
A_1(u^f,v^g)=2\int_{Y_{1T_1}\cap \{s\leq s_0\}} u_s^f\overline{v^g}dsdy'
\end{equation}
and $A_1(u^f,v^g)$ is uniquely determined by the D-to-N operator.
\qed

We shall construct a geometric optics solution of $L_1 u =0$
of the form (c.f. [E1]) :
\begin{equation}                        \label{eq:2.37}
u=u_N+u^{(N+1)},
\end{equation}
where
\begin{equation}                      \label{eq:2.38}
u_N=e^{ik(s-s_0)}\sum_{p=0}^N\frac{1}{(ik)^p}a_p(s,\tau,y').
\end{equation}
Substitute $u_N$ in (\ref{eq:2.17}) we get
\[
\frac{\partial a_0}{\partial \tau}-i\hat{A}_n^{(1)}a_0=0,\ \ 
4\frac{\partial a_p}{\partial \tau}-
4i\hat{A}_n^{(1)}a_p=-L_1a_{p-1},\ p\geq 1.
\]
We have
\[
a_0(s,\tau,y')=a_0(s,y')e^{ib(s,\tau,y')},
\]
where $b(s,\tau,y')=\int_{T-s}^\tau\hat{A}_n^{(1)}d\tau'.$
  We choose
$a_0(s,y')=\chi_1(s)\chi_2(y')$  where $\chi_1(s)\in C_0^\infty(\R^1),\ 
\chi_1(s)=1$ for $|s-s_0|<\delta,\ \chi_1(s)=0$ for $|s-s_0|>2\delta,
\ \chi_2(y')=\frac{1}{\e^{n-1}}\chi_0(\frac{y-y_0'}{\e}),\ \chi_0(y')\in
C_0^\infty(\R^{n-1}),\ \chi_0(y')=0$  for $|y'|>\delta,\ \ 
\int_{\R^{n-1}}\chi_0(y')dy'=1,\ \delta $ is small, $y_0'\in \Gamma$.
We define
\[
a_p=-\frac{1}{4}e^{ib}\int_{T-s}^\tau e^{-ib}L_1a_{p-1}d\tau',
\ \ 1\leq p\leq N,
\]
and we define $u^{(N+1)}$ as a solution of
\[
L_1u^{(N+1)}=-\frac{1}{4^N(ik)^N}(L_1a_N)e^{ik(s-s_0)},
\]
$u^{(N+1)}=u_t^{(N+1)}=0$  when  $t=T_1,\ u^{(N+1)}=0$  when
$y_n=0$.  Since $\mbox{supp\ }u_N$ is contained in a small neighborhood
of the line $\{s=s_0,y'=y_0'\}$ 
 we have that $u_N+u^{(N+1)}=u_{Nt}+u_t^{(N+1)}=0$  when 
$t=T_1, \ y_n >0$  and $\mbox{supp\ }(u_N+u^{(N+1)})\cap\{y_n=0\}\subset\Delta_{1T_1}.$

Substituting (\ref{eq:2.37}) in (\ref{eq:2.36}) we get that the principal
term in $k$ has the form
\[
ik\int_{Y_{1T_1}\cap\{s\leq s_0\} }e^{ik(s-s_0)}\chi_1(s)\chi_2(y')
e^{ib}\overline{v^g}dy'ds
\]
Note that 
$\tau=0$ on $Y_{1T_1}$.  Integrating by parts in $s$, taking the limit
when $k\rightarrow\infty$ and then taking the limit when
$\e\rightarrow 0$ we get that  $e^{ib(s_0,0,y_0')}\overline{v^g(s_0,0,y_0')}$
is determined by the D-to-N operator.
Here $(s_0,0,y_0')\in \Gamma\times [T_1,T]\subset Y_{1T_1}$ is 
arbitrary.  Replacing $T$ by $T',\ T'\in (T_1,T]$
we can determine $e^{ib}\overline{v^g(s,\tau,y')}$
for any $(s,\tau,y')\in X_0$   where $X_0=\{y'\in \Gamma,
T_1\leq s+\tau\leq T\}$.

Since we assumed that $v^g=0$ on $\partial Y_{1T_1}$ we get,  integrating by parts, 
that 
\[
A(u^f,v^g)=2\int_{Y_{1T_1}}u_s^f\overline{v^g}dsdy' 
= -2\int_{Y_{1T_1}}u^f\overline{v_s^g}dsdy'.
\]
Analogously to the proof of Lemmas \ref{lma:2.1} and \ref{lma:2.4}
one can prove that the integral
\begin{equation}                                \label{eq:2.39}
\int_{Y_{1T_1}\cap\{s\leq s_0\}}u^f\overline{v_s^g}dsdy' 
\end{equation}
is uniquely determined by the D-to-N operator.

Substitute the geometric optics solution (\ref{eq:2.37}) in (\ref{eq:2.39}).
Then integrating by parts in $k$, taking the limit when $k\rightarrow\infty$
and then the limit when $\e\rightarrow 0$ we get that $e^{ib}\overline{v_s^g}$
is determined by the D-to-N operator for any $(s,\tau,y')\in X_0$.  Since
we know $e^{ib}\overline{v^g}$ we know the derivative 
$\frac{\partial}{\partial s}(e^{ib}\overline{v^g})$.  We have 
\begin{equation}                               \label{eq:2.40}
\frac{\partial}{\partial s}(e^{ib}\overline{v^g})=
e^{ib}\overline{v_s^g}+ib_se^{ib}\overline{v^g}.
\end{equation}
Therefore we can find $b_s$ on the set $\{v^g\neq 0\}$.  
 Since, by Lemma \ref{lma:2.2} for $j=1$,
 $\{v^g\},\ g\in C_0^\infty(\Delta_{1T_1})$,  are dense in 
$\stackrel{o}{H^1}(R_{1T_1})$, where $R_{1T_1}=\Gamma\times[T_1,T]$,
 the union of all sets
   $\{v^g\neq 0\},\ g\in C_0^\infty(\Delta_{1T_1})$,  is dense in $R_{1T_1}$.
Therefore since $b_s $ is smooth
 we can 
recover $b_s$  on $R_{1T_1}$. 
Replacing $T$ by $T'\in (T_1,T]$ we can recover $b_s$  on $X_0$.
  Since $b=0$  when $y_n=0$ we can 
find  $b$ and consequently $v^g$.  Finally,  
$\hat{A}_n^{(1)}=\frac{\partial}{\partial \tau}b$. 

Therefore the D-to-N operator determines $\hat{A}_n^{(1)}(y,t)$ and $v^g$ in $X_0$
for all smooth
$g\in \stackrel{\circ}{H^1}(\Delta_{1T_1})$  where $v^g$ is the solution of
\begin{eqnarray}                                \label{eq:2.41}
L_1^*v^g=
\left(-i\frac{\partial}{\partial t} 
+\overline{\hat{A}_n^{(1)}}\right)^2v^g
-\left(-i\frac{\partial}{\partial y_n} 
+\overline{\hat{A}_n^{(1)}}\right)^2v^g
 \ \ \ \ \ \ \ \ \ \ \ \ \ \ \ 
\\
-\sum_{j,k=1}^{n-1}\left(-i\frac{\partial}{\partial y_j}
+\overline{\hat{A}_j^{(1)}}(y,t)
\right)
\hat{g}^{jk}(y)\left(-i\frac{\partial}{\partial y_k} +
\overline{\hat{A}_k^{(1)}}(y,t)\right)v^g
-\overline{\hat{V}}_1(y,t)v^g=0,
\nonumber
\end{eqnarray}
We can rewrite (\ref{eq:2.41}) in the form ( see (\ref{eq:2.17}) )
\begin{eqnarray}                               \label{eq:2.42}
-\sum_{j,k=1}^{n-1}\hat{g}^{jk}(y)v_{y_jy_k}^g(s,\tau,y')+
\sum_{j=1}^{n-1}B_j(s,\tau,y')v_{y_j}^g
\\
+C(s,\tau,y')v^g(s,\tau,y')=
4v_{st}^g-4i\overline{\hat{A}_n^{(1)}}v_s^g(s,\tau,y'),
\nonumber
\end{eqnarray}
where $B_j, C$ depend on $\hat{g}^{jk},\hat{A}_j^{(1)}, \hat{V}_1$.

Consider the restriction of (\ref{eq:2.42}) to $R_{1T_1}$,  i.e. when
$\tau=0$.
It follows from Lemma \ref{lma:2.2} that $\{v^g\},\ g\in 
\stackrel{\circ}{H^1}(\Delta_{1T_1}),$
are dense in $\stackrel{\circ}{H^1}(R_{1T_1})$.  
Pick arbitrary $v\in C_0^\infty(R_{1T_1})$.
Then there exists $v^{g_n}\in H^1(R_{1T_1})$ and smooth
such that $v^{g_n}\rightarrow v_0$
in $H^1(R_{1T_1})$.
Therefore $v^{g_n}\rightarrow v$ weakly in
$R_{1T_1}$, i.e.  $(v^{g_n},\varphi)\rightarrow(v,\varphi)$ for any 
$\varphi\in C_0^\infty(R_{1T_1})$.  Denote by $L_2v^g$ and $f^g$ the left hand side
and the right hand side of (\ref{eq:2.42}).
Since $v^{g_n}\rightarrow v$ weakly we get that  
$L_2v^{g_n}\rightarrow L_2v$ weakly.
Therefore $f^{g_n}=L_2v^{g_n}$ converges weakly to $f\stackrel{def}{=}L_2v$.
Note that $f^{g_n}$ is known since $v^{g_n}$ and $\hat{A}_n^{(1)}$ are known.  
Therefore we know $f$.   For any point $(s_0,y_0')$ we can find 
$M=\frac{n(n-1)}{2}+n-1+1=\frac{n(n+1)}{2}\ \ $ 
$C_0^\infty$ functions $v_1,...,v_M$  such that
the determinant $D(v_1,...,v_M)$ of the system (\ref{eq:2.42})  is not zero
at $(s_0,y_0')$.  Since $f_k=L_2v_k,k=1,...,M,$ are known at $(s_0,y_0')$
we can uniquely determine $\hat{g}^{jk}(y_0),B_j(s_0,0,y_0'),
C(s_0,0,y_0')$.
Since $(s_0,0,y'_0)$ is an arbitrary point of $R_{1T_1}$ we can find 
$\hat{g}^{jk}(y),\ \hat{A}_j^{(1)},\ 1\leq j\leq n-1,\ \hat{V}_1$ on $R_{1T_1}$.

Analoguosly, considering the intersection of $X_0$ with arbitrary plane $\tau=\tau_0$,
where $0\leq \tau_0\leq T-T_1$ we can determine the coefficients of
$L_1^*$ in $X_0$.

Therefore we proved the following theorem :
\begin{theorem}                         \label{theo:2.1}
The D-to-N operator on $\Gamma_0\times(0,T)$ uniquely determines the coefficients 
of $L_1$ in $X_0=\{y'\in \Gamma,T_1\leq s+\tau\leq T\}$.
\end{theorem}

\section{The global step.}
\label{section 3}
\init

The proofs in this section are similar to the proofs in  [E1].  Therefore
we will be brief.

Note that $D_0=\Gamma\times[0,\frac{T-T_1}{2}]$ is  the projection of $X_0$
on the plane $\{t=0\}$.   Since the coefficients of $L_1$ are 
analytic in $t$ the coefficients in $X_0$ determine the coefficients
of $L_1$ in $D_0\times(0,T_0)$.  We have the following result:

\begin{lemma}                        \label{lma:3.1}
The D-to-N operator on $\Gamma_0\times(0,T)$ determines the coefficients of 
$L_1$ in $D_0\times(0,T_0)$.
\end{lemma}

Let $L^{(p)},\ p=1,2,$ be two operators of the form (\ref{eq:1.1}) in domains 
$\Omega^{(p)}\times(0,T_0),\ p=1,2.$
Suppose $\Gamma_0\subset\partial\Omega^{(1)}\cap\partial\Omega^{(2)}$
and $\Lambda^{(1)}=\Lambda^{(2)}$ on $\Gamma_0\times(0,T_0)$  where
$\Lambda^{(p)}$ are the D-to-N operators corresponding to $L^{(p)},\ p=1,2$.

Let $y=\varphi_p(x)$ be the same as in (\ref{eq:2.1}),  $p=1,2.$
Denote $D^{(p)}=\varphi_p^{-1}(D_0)\subset\Omega^{(p)}$.
Then $\varphi=\varphi_1^{-1}\circ\varphi_2$  is a diffeomorphism of
$\overline{D^{(2)}}$ onto $\overline{D^{(1)}}$.
Let $L_1^{(p)}u_1^{(p)}=0$  in $D_0\times(0,T)$
where 
$L_1^{(p)}$ be the same as in (\ref{eq:2.16})  for  $p=1,2.$  Note that 
$u_1^{(p)}(y,t)$ are related to $u^{(p)}(x,t)$ where $L^{(p)}u^{(p)}=0$
by the formula (\ref{eq:2.19}) for $p=1,2$.  It follows 
from Remark 2.2 in [E1] that $\hat{g}_1(y',0)=\hat{g}_2(y',0)$.
Therefore there exists $\hat{c}(y,t)\in G_0(D_0\times(0,T_0))$,
analytic in $t$,
such that  $u^{(1)}(\varphi_1^{-1}(y),t)=\hat{c}(y,t)u^{(2)}(\varphi_2^{-1}(y),t)$ 
for $(y,t)\in \overline{D_0}\times(0,T_0)$.  Extend $\varphi(x)$ from
$\overline{D^{(2)}}$ to $\overline{\Omega^{(2)}}$ 
preserving the property
 that $\varphi=I$ on $\Gamma_0$ and $\varphi$ is a diffeomorphism 
of $\overline{\Omega^{(2)}}$ onto $\overline{\Omega_0}\stackrel{def}{=}
\overline{\varphi(\Omega^{(2)}}$.
The existence of such extension follows from [Hi],  Chapter 8. 
Let $c\in G_0(\Omega_0\times(0,T_0))$
be the extension analytic in $t$ of $\hat{c}(\varphi_2^{-1}(x),t)$ 
from $\overline{D^{(2)}}\times[0,T_0]$
to $\overline{\Omega_0}\times[0,T_0]$.
Then $L_0\stackrel{def}{=}c\circ \varphi\circ L^{(2)}$ is a differential operator
in $\Omega_0\times(0,T_0)$ such that $L_0=L^{(1)}$ in $D^{(1)}\times (0,T_0)$.
Let $B\subset D^{(1)}$ be a domain homeomorphic to a ball, 
$\overline{B}\cap\partial\Omega^{(1)}\stackrel{def}{=}S_1
\subset \Gamma$ and connected.  Let
$\Omega_1=\Omega_0\setminus\overline{B}$ and let 
$S_2=\partial B\setminus \overline{S_1}$.
    We assume that 
$\partial\Omega_1$ is smooth.   Denote by $\Lambda_0$ the D-to-N operator
corresponding to $L_0$.  Let $\delta=\max_{x\in B}d(x,\Gamma)$
where $d(x,\Gamma)$ is the distance in $\overline{B}$ from $x\in \overline{B}$ 
to $\Gamma$.

Consider $L_0,L^{(1)}$ and the corresponding D-to-N operators
$\Lambda_0,\Lambda^{(1)}$ in domains $(\Omega_0\setminus\overline{B})
\times(\delta,T_0-\delta),\ 
(\Omega^{(1)}\setminus\overline{B})\times(\delta,T_0-\delta)$,  respectively.

\begin{lemma}                              \label{lma:3.2}
If $\Lambda_0=\Lambda^{(1)}$ on $\Gamma_0\times(0,T_0)$ then
$\Lambda_0=\Lambda^{(1)}$  on 
$\Gamma_1\times(\delta,T_0-\delta)$ 
where $\Gamma_1=(\Gamma_0\setminus S_1)\cup S_2$.
\end{lemma}
Note that $\partial\Omega_1\setminus\gamma\subset\partial(\Omega^{(1)}\setminus
\overline{B})$
since $\Gamma_0\subset\partial\Omega_0\cap\partial\Omega^{(1)}$  and
$B\subset \Omega_0\cap\Omega^{(1)}$.
The proof of Lemma \ref{lma:3.2} is the same as the proof of Lemma 3.3
in [E1]  (c.f [KKL1],  Lemma 9).
The proof uses the Tataru's uniqueness theorem [T],  and this requires the 
analyticity of the coefficients in $t$.

Using repeatedly Lemmas \ref{lma:3.1} and \ref{lma:3.2} we get
a domain $\Omega^{(0)}\subset\Omega^{(1)},
\ \partial\Omega^{(0)}\cap\partial\Omega^{(1)} =\Gamma_0$,
a diffeomorphism $\varphi_2$ of $\overline{\Omega}^{(2)}$ onto
$\overline{\tilde{\Omega}}^{(2)}\stackrel{def}{=}
\varphi_2(\overline{\Omega}^{(2)})$
and a gauge transformation $\tilde{c}$ such that $\varphi_2=I$ on
$\Gamma_0,\ \tilde{\Omega}^{(2)}\supset\Omega^{(0)},\ \tilde{c}=1$
on $\Gamma_0\times(0,T_0)$ and $\tilde{c}\circ\varphi_2\circ L^{(2)}=L^{(1)}$
in $\Omega^{(0)}\times(0,T_0)$.
Moreover $\tilde{\Lambda}^{(2)}=\Lambda^{(1)}$ on $\partial\Omega^{(0)}\times
(T_0',T_0-T_0')$ for some $T_0'\in (0,\frac{T_*}{2})$
 where $\tilde{\Lambda}^{(2)}$ is the D-to-N operator 
corresponding  to $\tilde{L}^{(2)}\stackrel{def}{=}\tilde{c}\circ\varphi_2
\circ L^{(2)}$.

To complete the proof of Theorem \ref{theo:1.1}  we need two more lemmas.

Let $\gamma_0$ be an open subset of 
$\partial\Omega^{(0)}\setminus\overline{\Gamma_0}$
that is close to $\partial\Omega^{(1)}\setminus\overline{\Gamma_0}$ and
$\partial\tilde{\Omega}^{(2)}\setminus\overline{\Gamma_0}$
with respect to the corresponding metrics.  Let $\Delta_1$
be the union of all geodesics in $\overline{\Omega}^{(1)}\setminus\Omega^{(0)}$
starting at $\gamma_0$,  orthogonal to $\gamma_0$ and ending on 
$\partial\Omega^{(1)}$.
Analogously,  let $\tilde{\Delta}_2$  be the union of such geodesics in 
$\overline{\tilde{\Omega}}^{(2)}\setminus\Omega^{(0)}$.
We assume that the semi-geodesic coordinates can be introduced
in $\Delta_1$ and $\tilde{\Delta}_2$,   respectively.
Denote $\gamma_1=\overline{\Delta_1}\cap\partial\Omega^{(1)},\ \tilde{\gamma}_2= 
\overline{\tilde{\Delta}}_2\cap\partial\tilde{\Omega}^{(2)}$.

\begin{lemma}                              \label{lma:3.3}
There exists a diffeomorphism $\varphi^{(0)}$ of 
$\overline{\tilde{\Delta}}_2$ onto $\overline{\Delta_1}$ and 
a gauge $c^{(0)}$ on
$\overline{\Delta}_1\times (0,T_0)$  such that 
$\varphi^{(0)}(\tilde{\gamma}_2)=\gamma_1,\ \varphi^{(0)}=I$ on $\gamma_0,\ 
c^{(0)}=1$ on $\overline{\gamma}_0\times(0,T_0)$ and 
$c^{(0)}\circ\varphi^{(0)}\circ\tilde{L}^{(2)}=L^{(1)}$ in $\overline{\Delta_1}\times
(0,T_0)$ .
\end{lemma}
 
\qed

Consider  a situation
when there exists a part $\Gamma_1$ of $\partial\Omega^{(0)}$ that is close 
to a different part of $\partial\Omega^{(0)}$.  Denote by $D_1$ the 
union of all geodesics in $\overline{\Omega^{(1)}}\setminus\Omega^{(0)}$
starting on $\Gamma_1$,  orthogonal to $\Gamma_1$ and ending on
a set $\Gamma_2\subset\partial\Omega^{(0)}$.  Analogously let
$D_2$ be the union of the geodesics in 
$\overline{\tilde{\Omega}^{(2)}}\setminus\Omega^{(0)}$,
starting at $\Gamma_1$,  orthogonal to $\Gamma_1$ and ending on
some set $\Gamma_2'\subset\partial\Omega^{(0)}$.

As before we assume that the semi-geodesic coordinates are introduced
in $D_1$ and $D_2$.

\begin{lemma}                                  \label{lma:3.4}
There exists a diffeomorphism $\psi$ of $\overline{D_2}$ onto 
$\overline{D_1}$ and a gauge $c$ on $\overline{D_1}\times(0,T_0)$
such that $\psi=I$ on $\Gamma_1,\ \psi=I$ on $\Gamma_2$,  in particular,
$\Gamma_2'=\Gamma_2,\ c=1$  on $\Gamma_1\times(0,T_0),$ and  on
$\Gamma_2\times(0,T_0)$.  Moreover,
$c\circ\psi\circ L^{(2)}=L^{(1)}$  in $D_1\times(0,T_0)$ .
\end{lemma}

The proofs of Lemmas  \ref{lma:3.3} and \ref{lma:3.4} 
are the same as of corresponding results in  [E1].
Combining Lemmas \ref{lma:3.1}, \ref{lma:3.2}, \ref{lma:3.3}, \ref{lma:3.4}
we prove Theorem \ref{theo:1.1} (c.f. [E1]).

\section{Proof of Lemma \ref{lma:2.2}.}
\label{section 4}
\init

Denote by $\Delta_1$ a domain in $\R^{n+1}$ bounded by three planes :
$\Gamma_2=\{\tau=T-t-y_n=0,\ 0\leq y_n\leq \frac{T-T_1}{2}\},
\ \Gamma_3=\{s=t-y_n=T_1,\ \frac{T+T_1}{2}\leq t\leq T\},
\ \Gamma_4=\{t=T,\ 0\leq y_n\leq T-T_1\}$. 
Let $\mathcal{H}=H^1(\Gamma_4)\times L_2(\Gamma_4)$
and
let $\mathcal{H}_1$ be the space of pairs $\{\varphi,\psi\},\ 
\varphi\in H^1(\Gamma_2),\ \psi=H^1(\Gamma_3),
\ \varphi=\psi$ when $t=\frac{T+T_1}{2}$, with
the norm $\|\{\varphi,\psi\}\|_1^2=\|\varphi\|_{1,\Gamma_2}^2+
\|\psi\|_{1,\Gamma_3}^2.$  We shall assume that all functions have 
a compact support in $y'\in \R^{n-1}$.

The following lemma is an extension of Lemma 3.1  in [E1] :
\begin{lemma}                                \label{lma:4.1}
For any $\{v_0,v_1\}\in \mathcal{H}_1$ there exists $\{w_0,w_1\}\in
\mathcal{H}$ and $u\in H^1(\Delta_1)$ such that $L_1u=0$ in
$\Delta_1,\ u\left|_{\Gamma_2}=v_0\right., \ u\left|_{\Gamma_3}=v_1\right.,
\ u\left|_{\Gamma_4}=w_0\right.,
\ \frac{\partial u}{\partial t}\left|_{\Gamma_4}=w_1\right.$.
And vice versa,  for any $\{w_0,w_1\}\in \mathcal{H}$
there exists $\{v_0,v_1\}\in \mathcal{H}_1$ and $u\in H^1(\Delta_1)$ with
the same properties.  The iinequalities (\ref{eq:4.2}) and (\ref{eq:4.3})
hold. 
\end{lemma}

{\bf Proof:} Let $\Delta_{1,T'}$ be the domain bounded by $\Gamma_2,
\Gamma_3$ and $\Gamma_{4,T'}$  where $\Gamma_{4,T'}$ is the plane 
$\{t=T'\},\ T'\in[\frac{T+T_1}{2},T]$.  Let $\|u\|_{1,T'}^2=
\|u_t\|_{0,\Gamma_{4,T'}}^2+\|u\|_{1,\Gamma_{4,T'}}^2$.
Denote by $(u,v)_{\Delta_{1,T'}}$  the $L_2$-inner product in 
$\Delta_{1,T'}$.
Let $u$ be smooth in $\overline{\Delta_{1,T}}$, has compact support
in $y'$ and $L_1u=0$ in $\Delta_{1,T}$.  Integrating by parts 
$0=(L_1u,u_t)_{\Delta_{1,T'}}+(u_t,L_1u)_{\Delta_{1,T'}}$
over $\Delta_{1,T'}$  (c.f. (3.1), (3.2) in [E1])   we get
\begin{equation}                               \label{eq:4.1}
\|u\|_{1,T'}^2\leq C(\|u|_{\Gamma_{2,T'}}\|_{1,\Gamma_{2,T'}}^2+
\|u|_{\Gamma_{3,T'}}\|_{1,\Gamma_{3,T'}}^2
+\int_{\frac{T+T_1}{2}}^{T'}\|u\|_{1,t}^2dt),
\end{equation}
where $\Gamma_{2,T'},\Gamma_{3,T'}$ are parts of $\Gamma_2,\Gamma_3$
respectively where $t\leq T'$.  Since  $T-T_1$ is small we get
\begin{equation}                               \label{eq:4.2}
\max_{T'\in [\frac{T}{2},T]}\|u\|_{1,T'}^2\leq C_1(\|u\|_{1,\Gamma_{2,T}}^2+
\|u\|_{2,\Gamma_{3,T}}^2),
\end{equation}
where $\Gamma_{2,T}=\Gamma_2,\ \Gamma_{3,T}=\Gamma_3$.
In particular, 
\[
\|u\|_{1,\Gamma_4}^2+\|u_t\|_{0,\Gamma_4}^2
\leq C(\|u\|_{1,\Gamma_2}^2+\|u\|_{1,\Gamma_3}^2).
\]
The reverse inequality also holds:
\begin{equation}                         \label{eq:4.3}
\|u\|_{1,\Gamma_2}^2+\|u\|_{1,\Gamma_3}^2
\leq C(\|u\|_{1,\Gamma_4}^2+\|u_t\|_{0,\Gamma_4}^2),
\end{equation}
where $L_1u=0$ in $\Delta_1,\ u$ is smooth.

To prove this denote by $\Gamma_{2,\tau'}$ the intersection of
$\overline{\Delta_1}$ with the plane $\tau=\tau',\ 0\leq -\tau'\leq T-T_1$ and 
denote by $\Gamma_{3,s'}$ the intersection of $\overline{\Delta_1}$ with 
the plane $s=s',\ T_1\leq s'\leq T$.  Integrating by parts 
$0=(L_1u,u_t)_{\Delta_{s'}}+(u_t,L_1u)_{\Delta_{s'}}$ we get, as in 
(\ref{eq:4.1}) (c.f. [E1]),
\begin{equation}                        \label{eq:4.4}
\|u\|_{1,\Gamma_{3,s'}}^2 +\|u\|_{1,\Gamma_{2,s'}}^2
\leq C(\|u\|_{1,\Gamma_4}^2 + \|u_t\|_{0,\Gamma_4}^2+I_1),
\end{equation}
where 
\[
I_1=\int_{\Delta_{s'}}\left(|u_s|^2+|u_\tau|^2+
\sum_{j=1}^{n-1}|u_{y_j}|^2+|u|^2\right)dy'dsd\tau.
\]
Here
$\Delta_{s'}$ is the part of $\Delta_1$ where  $s>s'$  and 
$\Gamma_{2,s'}$ is the part of $\Gamma_2$ where $s>s'$.  Let
$\Delta_{\tau'}$ be the part of $\Delta_1$  where $\tau<\tau'$.  Then
integrating by parts $0=(L_1u,u_t)_{\Delta_{\tau'}} + 
(u_t,L_1u)_{\Delta_{\tau'}}$ we get 
\begin{equation}                        \label{eq:4.5}
\|u\|_{1,\Gamma_{2,\tau'}}^2 +\|u\|_{1,\Gamma_{3,\tau'}}^2
\leq C(\|u\|_{1,\Gamma_4}^2 + \|u_t\|_{0,\Gamma_4}^2+I_2),
\end{equation}
where 
$I_2$ is the same as $I_1$ with $\Delta_{s'}$ replaced by
$\Delta_{\tau'},\ \Gamma_{3,\tau'}$  is the part of 
$\Gamma_3,$ where $\tau <\tau'$.
Since
\[
I_1+I_2\leq C(T-T_1)(\max_{s'}\|u\|^2_{1,\Gamma_{3,s'}} +
\max_{\tau'}\|u\|^2_{1,\Gamma_{2,\tau'}},
\]
we get from (\ref{eq:4.4}), (\ref{eq:4.5})
\[
\max_{s'}\|u\|^2_{1,\Gamma_{3,s'}}+\max_{\tau'}\|u\|^2_{1,\Gamma_{2,\tau'}}
\]
\[
\leq
C_1(T-T_1)(\max_{s'}\|u\|^2_{1,\Gamma_{3,s'}}+\max_{\tau'}\|u\|^2_{1,\Gamma_{2,\tau'}})
+C(\|u\|_{1,\Gamma_4}^2+\|u_t\|_{0,\Gamma_4}^2).
\]
Since $T-T_1$ is small we get (\ref{eq:4.3}).

Let $\{\varphi_0,\psi_0\}$ be a smooth pair belonging to
$\mathcal{H}_1$.
Define $b(y',y_n,t)=\varphi_0(s,y')+\psi_0(\tau,y')-\varphi_0(0,y').$
Note that 
$\varphi_0(0,y')=\psi_0(0,y'),\ \ b|_{\Gamma_2}=\varphi_0(s,y'),
b|_{\Gamma_3}=\psi_0(\tau,y')$.

Since $b(y',y_n,t)$ is smooth in $\overline{\Delta_1}$ we have that $L_1b$
is smooth in $\overline{\Delta_1}$.  Let $f=-L_1b\ \ $   in $\overline{\Delta_1},
\ f=0$ otherwise.  Then $f=0$ for $t<\frac{T-T_1}{2},
\ \int_{\frac{T+T_1}{2}}^T\|f\|_{0,\Gamma_{4,t}}dt<+\infty$.
Let $u_0$ be the solution of $L_1u_0=f$  when 
$t>\frac{T+T_1}{2}$ with zero initial conditions when $t=\frac{T+T_1}{2}$.
It is well-known (see, for example, [H])  that there exists a unique such $u_0$
and 
\begin{equation}                           \label{eq:4.6}
\|u_0\|_{1,T}\leq C\int_{\frac{T+T_1}{2}}^T\|f\|_{0,\Gamma_{4,t}}dt.
\end{equation}

Since $u_0$ has zero initial conditions when $t=\frac{T+T_1}{2}$ and $f=0$
outside of $\Delta_1$ we get by the domain of dependence argument that
$u_0=0$ outside $\overline{\Delta_1},\ t\leq T$.  Since the restrictions
of $u_0$ on the planes $s=\mbox{const}$
and $\tau=\mbox{const},\ t\leq T$, are continuous in $s$ and $\tau$ 
respectively we get
that 
\[
u_0|_{\Gamma_2}=0,\ \ u_0|_{\Gamma_3}=0.
\]
Denote $u=b+u_0$.  Then
\begin{equation}                              \label{eq:4.7}
L_1u=0\ \ \ \mbox{in}\ \ \Delta_1,
\end{equation}
\begin{equation}                           \label{eq:4.8}
u|_{\Gamma_2}=\varphi_0,\ \ u|_{\Gamma_3}=\psi_0
\end{equation}
and
\begin{equation}                          \label{eq:4.9}
\{u,u_t\}|_{\Gamma_4}\in \mathcal{H}.
\end{equation}

Let $\{v_0,v_1\}\in \mathcal{H}_1$ be arbitrary.  Take a sequence
$\{\varphi_n,\psi_n\}\in\mathcal{H}_1$ of smooth functions such that
$\{\varphi_n,\psi_n\}\rightarrow \{v_0,v_1\}$  in $\mathcal{H}_1$.
Let $u_n$ be a sequence such that (\ref{eq:4.7}), (\ref{eq:4.8}), (\ref{eq:4.9})
hold with $\varphi_0,\psi_0$ replaced by $\varphi_n,\psi_n$.  Then (\ref{eq:4.2})
implies that $u_n$ converges to $u$ in $\max_{\frac{T+T_1}{2}\leq 
T'\leq T} \|u\|_{1,T'}$
norm and $L_1u=0$ in $\Delta_1,\ u|_{\Gamma_2}=v_0,\  u|_{\Gamma_3}=v_1,
\ \{w_0,w_1\}\in \mathcal{H}$  where $w_0=u|_{\Gamma_4},\ w_1=u_t|_{\Gamma_4}$.

Therefore the map $\{v_0,v_1\}\rightarrow\{w_0,w_1\}$ is a bounded map of
$\mathcal{H}_1$ to $\mathcal{H}$.  Take any smooth pair 
$\{\varphi^{(0)},\psi^{(0)}\}\in\mathcal{H}$.  Solving the Cauchy problem
with the inithial data $\{\varphi^{(0)},\psi^{(0)}\}$ (see, for example, [H])
we get a smooth pair in $\mathcal{H}_1$.  Therefore the image of the map
$\mathcal{H}_1\rightarrow \mathcal{H}$ is dense.  Using the estimate 
(\ref{eq:4.3}) we get that the map $\mathcal{H}_1\rightarrow \mathcal{H}$ 
in one-to-one and onto.
\qed

Denote by $\Delta_2$ the domain bounded by the plane 
$\Gamma_2'=\{\tau=T-y_n-t=0,\ T_1< t< T\}$, by the plane $\{y_n=0,\ T_1<t<T\}$ and
by the plane $\{t=T_1,\ 0<y_n< T-T_1\}$. Let $\Gamma_\infty$ be 
the plane $\{\tau =0\}$.   Denote by
$\stackrel{\circ\ \  }{H^{-1}}(\Gamma_2')\subset H^{-1}(\Gamma_\infty)$  
the Sobolev space of
distributions in $\Gamma_\infty$ having supports in $\overline{\Gamma_2'}$
and by $H^{-1}(\Gamma_2')$  the space of the restrictions of 
distributions from $H^{-1}(\Gamma_\infty)$ to $\Gamma_2'$.
\begin{lemma}                                \label{lma:4.2}
Let $h\in H^{-1}(\Gamma_2')$ be arbitrary.
  There exists a distribution $u$ in $\Delta_2$ such that  
$L_1u=0$ in $\Delta_2,\ u$ has a restriction to $y_n=0,
\ \frac{\partial u}{\partial s}$ has a restriction to $\Gamma_2'$ and
$u|_{y_n=0}=0,\ \frac{\partial u}{\partial s}|_{\Gamma_2'}=h$.
\end{lemma}

{\bf Proof:}  Let $h_0\in H^{-1}(\Gamma_\infty)$ be an extension of $h$
to $\Gamma_\infty$.
We always can choose $h_0$ such 
that $\mbox{supp\ } h_0\subset \overline{\Gamma_2'}$,  
i.e. $h_0\in \stackrel{\circ\ \ \ }{H^{-1}}(\Gamma_2')$.
Let $v(s,y')$ be such that $v=0$ for $s>T-T_1$ and
$\frac{\partial v}{\partial s}=h_0$.  We have that $v\in H^{0,-1}$,
i.e. $v$ belongs to $L_2$ in $s$ and $v$ belongs to $H^{-1}$ in $y'$.

Denote by $\Delta_2^-$ the reflection of $\Delta_2$
with respect to the plane $y_n=0$ and let $\Delta_3=
\Delta_2\cup\Delta_2^-\cup \{y_n=0, T_1<t <T\}$.
Extend  for $y_n< 0$ the coefficient before
$\frac{\partial}{\partial y_n}$ in (\ref{eq:2.16})
as an odd function in
$\Delta_3$  and extend the remaining coefficients of $L_1$ as even 
functions.  We shall construct an odd distribution solution of $L_1u=0$
in $\Delta_3$.  Then automatically $u|_{y_n=0}=0$.  Note that any distribution
solution of $L_1u=0$ in $\Delta_3$ is continuously differentiable in
$y_n$ as a function of $y_n$ with the values that are distributions
in $(y',t)$.  This is a consequence of the fact that $y_n=0$ is not
a characteristic plane.  Therefore the restrictions 
$u|_{y_n=0}$ and $\frac{\partial u}{\partial y_n}|_{y_n=0}$ exist.

Denote by $\hat{\Gamma}_2'$ the reflection of $\Gamma_2'$ with respect to
$y_n=0$ and let $v_1(\tau,y')$ on $\hat{\Gamma}_2'$ 
be such that $v(s,y')+v_1(\tau,y')=0$
when $y_n=0$,  i.e. $v+v_1$ is odd with respect to $y_n$ in $\Delta_3$.
We shall look for the solution $L_1u=0$ in the form
\[
u=v(s,y')+v_1(\tau,y')+w.
\]
Then $w$ satisfies the equation $L_1w=g$,  where $g=-L_1(v+v_1)$.  Note
that $L_1$ has the form (\ref{eq:2.17})  for $y_n>0$.  Therefore
$g=g_1+g_2$,  where $g_1=4i\hat{A}_nv_s$ for $y_n>0$, $g_1$ is
odd in $y_n$,  and  $g_2=-L_1'(v+v_1)$ where $L_1'$ is a differential operator 
in $y'$ of order 2.  Since $v$ and $v_1$ belong to $L_2$ in $s$ and $\tau$
and to $H^{-1}$ in $y'$ we have that $g_2\in H^{0,-3}$,  i.e. $g_2$
belongs to $L_2$ in $s$ and $\tau$ and to $H^{-3}$ in $y'$.
Denote by $g_{20}$ the extension  of $g_2$ by zero outside of $\Delta_3$
for $t>T_1$.  Let $w_2$ be the solution of the Cauchy problem 
$L_1w_2=g_{20}$ for $t>T_1,\ w_2=0$ for $t>T$.  Since $g_{20}\in H^{0,-3}$
there exists a unique solution $w_2\in H^{1,-3}$  (see, for example, [E4] or [H]).
Note that $w_2$ is odd in $y_n$ since $g_{20}$ is odd.  Since $w_2=0$ for
$t>T$ we get that $w_2|_{\Gamma_2'}=0,\ w_2|_{\hat{\Gamma}_2'}=0$ by the domain of
dependence argument.

Denote by $g_{10}$ 
the odd extension in $y_n$ of $g_1\theta(\tau)$ given
for $y_n>0$ where $\theta(\tau)=1$ for $\tau>0,\ 
\theta(\tau)=0$ for $\tau<0$.  Then $g_{10}=0$ outside of $\Delta_3$ 
for $t>T_1$  since $v(s,y')=0$ for $s>T-T_1$.  Let $w_1$ be the unique odd
solution of the Cauchy problem $L_1w_1=g_{10}$ for $t>T_1,\ w_1=0$ for
$t>T$ (see [H] or [E4]).

By the domain of dependence argument and 
since $w_1$ is continuous in $\tau$ near $\Gamma_2'$ we get
that $w_1|_{\Gamma_2'}=0$.  Therefore $\frac{\partial w_1}{\partial s}|_{\Gamma_2'}=0$.
Therefore $u=v(s,y')+v_1(\tau,y')+w_1+w_2$ satisfies $L_1u=0$ in
$\Delta_2,\ u|_{y_n=0}=0$ and $\frac{\partial u}{\partial s}=
\frac{\partial v}{\partial s}=h$ on $\Gamma_2'$ since 
$\frac{\partial w_1}{\partial s}=\frac{\partial w_2}{\partial s}=
\frac{\partial v_1}{\partial s}=0$ on $\Gamma_2'$.
\qed

Now we are ready to conclude the proof of Lemma \ref{lma:2.2}.
As in [E1] ( see the beginning of  \S 3 in [E1]),  in order to prove
that $\{v^g\},\ g\in H_0^1(\Delta_{js_0})$, 
are dense in $H_0^1(R_{js_0}),j=1,2,3,$ it is enough
to prove that
$\{v^g\},\ g\in C_0^\infty(\Delta_{js_0})$, are dense in 
$\stackrel{\circ}{H^1}(R_{js_0})$.  
Fix,  for definiteness,  $j=3$.  Suppose that there exists 
$v_1\in \stackrel{\circ}{H^1}(R_{3s_0})$ and 
 $h\in H^{-1}(Y_{3s_0})$  such
that $(h,v^g)=0$ for all $g\in C_0^\infty(\Delta_{3s_0})$ and
$(h,v_1)=1$. Note that $(h,v)$,  where $v\in \stackrel{\circ\ \ }{H^1}(\Gamma_2')$,
is understood as  $(lh,v)$,  where $lh\in H^{-1}(\Gamma_\infty)$ is an
arbitrary extension of $h$.

Let $u$ be the same as in Lemma \ref{lma:4.2}.  Applying Green formula
(\ref{eq:2.24}) to $u$ and $v^g$ we get
\begin{equation}                         \label{eq:4.10}
0=2(\frac{\partial u}{\partial s},v^g)=\int_{\Delta_{3s_0}}
\frac{\partial u_1}{\partial y_n} \overline{g}dy'dt.
\end{equation}
To justify the Green formula (\ref{eq:4.10}) one can take a sequence 
$h_n\in C_0^\infty(\Gamma_2')$ such that $h_n\rightarrow h_0$ in
$\stackrel{\circ\ \ \ }{H^{-1}}(\Gamma_2')$,  where $h_0$ is
an extension of $h$.  By Lemma \ref{lma:4.2}
there exists $u_n$ such that $L_1u_n=0$ in $\Delta_2,\ u_n|_{y_n=0}=0,\ 
 \frac{\partial u_n}{\partial s}|_{\Gamma_2'}=h_n$.
Note that $(u_{ns},v^g)=-(u_n,\frac{\partial v^g}{\partial s})$ since
$v^g=0$ on $\partial Y_{1s_0}$ and
\begin{equation}                      \label{eq:4.11}
2(u_{ns},v^g)
=\int_{\Delta_{3s_0}}\frac{\partial u_n}{\partial y_n}\overline{g}
dy'dt
\end{equation}

Passing to the limit in (\ref{eq:4.11}) we get (\ref{eq:4.10}).  Since
$g\in C_0^\infty(\Delta_{3s_0})$ is arbitrary we get that 
$\frac{\partial u}{\partial y_n}=0$ in $\Delta_{3s_0}$.  Extend $u$
by zero for $y_n<0,\ s_0<t<T$.  By the Tataru's uniqueness theorem [T] (note
that the theorem holds when $u$ is a distribution too [T1])  we get 
$u=0$ in  the double cone
of influence of $\Delta_{3s_0}$. 
 In particular,  $h=\frac{\partial u}{\partial s}=0$
on the interior of $R_{3s_0}$. 
Therefore $(h,v_1)=0$  since $\mbox{supp\ } v_1\subset \overline{R}_{3s_0}$ and
this contradicts the assumption that $(h,v_1)=1$. 
   The proof for $j=1,2$ is 
identical.
\qed


\begin{thebibliography}{9999}
\bibitem[BLR]{} Bardos, C., Lebeau, G. and Rauch, J., 1992,
Sharp sufficient conditions for the observation, control and
stabilization of waves from the boundary,
SIAM J. Contr. Opt. 30, 1024-1065
\bibitem[B1]{} Belishev, M., 1997, Boundary control in reconstruction
of manifolds and metrics (the BC method),
Inverse Problems 13, R1-R45
\bibitem[B2]{} Belishev, M., 2002, How to see waves under the Earthsurface
 (the BC-method for geophysicists), Ill -Posed and Inverse Problems,
55-72 (S.Kabanikhin and V.Romanov (Eds), VSP)
\bibitem[B3]{} Belishev, M., 1997, On the uniqueness of the 
reconstruction of lower-order terms of the wave equation from dynamic
boundary data, (Russian), 
Zapiski Nauchnih Seminarov POMI, 29,  55-76
\bibitem[E1]{} Eskin, G., 2005, A new approach to the hyperbolic
inverse problems, ArXiv:math.AP/0505452 
\bibitem[E2]{} Eskin, G., 2005, Inverse problems for Schr\"{o}dinger
equations with Yang-Mills potentials in domains with obstacles and
the Aharonov-Bohm effect, Institute of Physics Conference Series 12,
23-32, ArXiv:math.AP/0505554
\bibitem[E3]{} Eskin, G., 2004,  Inverse boundary value problems in
domains with several obstacles, Inverse Problems 20, 1497-1516
\bibitem[E4]{} Eskin, G., 1987,  Mixed initial-boundary value problems 
for second order hyperbolic equations, Comm. in PDE, 12, 503-87
\bibitem[H]{} Hormander, L., 1985, The Analysis of Linear Partial 
Differential Operators III (Berlin: Springer)
\bibitem[Hi]{} Hirsch, M., 1976, Differential Topology (New York:Springer)
\bibitem[I]{} Isakov, V., 1998, Inverse problems for partial differential 
equations,  Appl. Math. Studies, vol. 127, Springer, 284 pp.
\bibitem[KKL]{} Katchalov, A., Kurylev, Y., Lassas, M., 2001,
Inverse boundary spectral problems (Boca Baton : Chapman\&Hall)
\bibitem[KKL1]{} Katchalov, A., Kurylev, Y., Lassas, M., 2004,
Energy measurements and equivalence of boundary data for inverse problems
on noncompact manifolds,  IMA Volumes, v.137, 183-214
\bibitem[KK]{} Katchalov, A., Kurylev, Y., 1998,
Multidimensional inverse problems with incomplete boundary spectral data,
Comm. Part. Diff. Eq. 23, 55-95
\bibitem[K]{} Kurylev, Y., 1993, Multi-dimensional inverse boundary problems
by BC-mathod : groups of transformations and uniqueness results,
Math. Comput. Modelling 18, 33-45
\bibitem[KL1]{} Kurylev, Y. and Lassas, M., 2000,
Hyperbolic inverse problems with data on a part of the boundary
AMS/1P Stud. Adv. Math, 16, 259-272
\bibitem[KL2]{} Kurylev, Y. and Lassas, M., 2002,
Hyperbolic inverse boundary value problems and time-continuation 
of the non-stationary Dirichlet-to-Neumann map,
Proc. Royal Soc. Edinburgh, 132, 931-949
\bibitem[KL2]{} Kurylev, Y. and Lassas, M.,1997,
The multidimensional Gel'fand inverse problem for nonself-adjoint
operators,  Inverse Problems, 13, 1495-1501
\bibitem[RS]{} Ramm, A. and Sjostrand, J., 1991, An inverse problem
of the wave equation, Math. Z., 206, 119-130
\bibitem[St]{} Stefanov, P., 1989, Uniqueness of multidimensional 
inverse scattering problem with time-dependent potentials,
Math. Z., 201,  541-549
\bibitem[T]{} Tataru, D., 1995, Unique continuation for solutions
to PDE, Comm. in PDE 20, 855-84
\bibitem[T1]{} Tataru, D., Private communication.



\end{thebibliography}
\end{document}